\documentclass[10pt]{amsart}
\allowdisplaybreaks[4]
\usepackage{color}

\definecolor{c20}{rgb}{0.,0.7,0.}
\definecolor{c30}{rgb}{0.,0.,1.}
\definecolor{c40}{rgb}{1,0.1,0.7}
\definecolor{c50}{rgb}{1,0,0}
\definecolor{c60}{rgb}{1,0.9,0.1}

\def\Pe#1{\textcolor{c20}{#1}}
\def\Pe#1{#1}

\def\Ke#1{\textcolor{c20}{#1}}
\def\Ke#1{#1}

\def\De#1{\textcolor{c20}{#1}}
\def\De#1{#1}

\def\Ee#1{\textcolor{c20}{#1}}
\def\Ee#1{#1}
\def\Hh#1{\textcolor{c20}{#1}}
\def\Hh#1{#1}

\def\cLP#1{\textcolor{c50}{#1}}
\def\cLP#1{#1}
\def\eLP#1{\textcolor{c50}{#1}}
\def\eLP#1{#1}

\def\eP#1{\textcolor{c50}{#1}}
\def\eP#1{#1}

\newcommand{\abs}[1]{\lvert #1 \rvert}

\newcommand{\E}[1]{\mathbb{E}\left\{#1\right\}}

\newcommand{\pk}[1]{\mathbb{P} \left\{ #1 \right\} }

\newcommand{\R}{\mathbb{R}}

\newcommand{\N}{\mathbb{N}}
\newcommand{\inr}{\in \R}
\newcommand{\inn}{\in \N}

\newcommand{\limit}[1]{\lim_{#1 \to   \infty}}

\newcommand{\BQN}{\begin{eqnarray}}
\newcommand{\EQN}{\end{eqnarray}}
\newcommand{\BQNY}{\begin{eqnarray*}}
\newcommand{\EQNY}{\end{eqnarray*}}

\newcommand{\BS}{\begin{sat}}
\newcommand{\ES}{\end{sat}}
\newcommand{\BT}{\begin{theo}}
\newcommand{\ET}{\end{theo}}
\newcommand{\BK}{\begin{korr}}
\newcommand{\EK}{\end{korr}}

\newcommand{\BD}{\begin{de}}
\newcommand{\ED}{\end{de}}

\newcommand{\BIT}{\begin{itemize}}
\newcommand{\EIT}{\end{itemize}}
\newcommand{\BDI}{\begin{description}}
\newcommand{\EDI}{\end{description}}

\newcommand{\BRM}{\begin{remarks}}
\newcommand{\ERM}{\end{remarks}}

\newcommand{\BEL}{\begin{lem}}
\newcommand{\EEL}{\end{lem}}

\newtheorem{theo}{Theorem}[section]
\newtheorem{sat}[theo]{Proposition}
\newtheorem{de}[theo]{Definition}
\newtheorem{lem}[theo]{Lemma}

\newtheorem{korr}[theo]{Corollary}
\newtheorem{remark}[theo]{Remark}
\newtheorem{remarks}[theo]{Remarks}

\newcommand{\netheo}[1]{{Theorem \ref{#1}}}

\newcommand{\prooftheo}[1]{ \textsc{\bf Proof of Theorem} \ref{#1}:}

\newcommand{\prooflem}[1]{\textsc{\bf Proof of Lemma} \ref{#1}:}

\newcommand{\COM}[1]{}

\newcommand{\QED}{\hfill $\Box$}

\def\rw{\rightarrow}

\def\IF{\infty}
\def\piter{\Ee{\mathcal{Q}}}

\def\oo{(1+o(1))}
\def\asu{\ \ \text{as}\ u\rw\IF}
\def\LT{\left}
\def\RT{\right}
\def\H{\mathcal{H}}

\def\ooo{(1+o(1))}

\def\rw{\rightarrow}

\def\Del{\triangle}
\def\Delo{\triangle^1}
\def\Delt{\triangle^2}

\def\vn{\varepsilon}
\def\Var{\mathbb{V}ar}

\def\hs{\hat{s}_0}
\def\htt{\hat{t}_0}
\def\Cov{\mathbb{C}ov}
\def\hatc{\hat{c}}

\def\tAu{\tau_1^*(u)}
\def\tAuu{\tau_2^*(u)}
\newcommand{\equaldis}{\stackrel{d}{=}}
\newcommand{\todis}{\stackrel{d}{\to}}
\def\LamT{\De{ \Lambda_T}}

\begin{document}

\title[Extremes and First Passage Times of Correlated  fBm's] {Extremes and First Passage Times of Correlated\\ Fractional Brownian Motions}

\author{Enkelejd  Hashorva}
\address{Enkelejd Hashorva, University of Lausanne\\
B\^{a}timent Extranef, UNIL-Dorigny, 1015 Lausanne, Switzerland
}
\email{Enkelejd.Hashorva@unil.ch}

\author{Lanpeng Ji }
\address{Lanpeng Ji, University of Lausanne\\
B\^{a}timent Extranef, UNIL-Dorigny, 1015 Lausanne, Switzerland
}
\email{Lanpeng.Ji@unil.ch}

\bigskip

 \maketitle

{\bf Abstract:} Let $\{X_i(t),t\ge0\}, i=1,2$ be two \De{standard} fractional Brownian motions \Pe{being jointly Gaussian} with constant cross-correlation.
In this paper we derive the exact asymptotics of the \De{joint survival function}
$$
\pk{\sup_{s\in[0,1]}X_1(s)>u,\ \sup_{t\in[0,1]}X_2(t)>u}
$$
as $u\rightarrow \infty$. \Ke{A novel finding of this contribution is the exponential approximation of the
joint conditional first passage times of $X_1, X_2$}. As a by-product we obtain generalizations of the
Borell-TIS inequality and the Piterbarg inequality for 2-dimensional Gaussian random fields.

{\bf Key Words:} Extremes; first passage times;
Borell-TIS inequality; Piterbarg inequality; fractional Brownian motion; Gaussian random fields.

{\bf AMS Classification:} Primary 60G15; secondary 60G70

\section{Introduction and main \eLP{results}}

Let $\{X_i(t),t\ge0\}, i=1,2$ be two standard fractional Brownian motion's (fBm's) with \eP{Hurst} indexes $\alpha_i/2\in(0,1),i=1,2,$ i.e., \Hh{$X_i$ is a}  \Ee{centered} Gaussian process with \Hh{a.s.\ continuous sample paths and}
covariance  function
\BQNY
\mathbb{C}ov(X_i(t),X_i(s))=\frac{1}{2}(\abs{t}^{\alpha_i}+\abs{s}^{\alpha_i}-\mid t-s\mid^{\alpha_i}),\quad s,t\ge0,\ \ \ \ i=1,2.
\EQNY
\Pe{Hereafter $(X_1,X_2)$ are assumed to be jointly Gaussian with cross-correlation function} 
$r(s,t)\eLP{=\E{X_1(s)X_2(t)}/\sqrt{s^{\alpha_1}t^{\alpha_2}}}\in(-1,1)$. \Pe{Calculation of the following joint survival function}
\BQN\label{eq:Pru}
P_r(u):=\pk{\sup_{s\in[0,1]}X_1(s)>u,\ \sup_{t\in[0,1]}X_2(t)>u},\ \ \ \ u>0.
\EQN
 is important for various applications in   statistics,
mathematical finance and insurance mathematics.  
\De{The special} simple model of two correlated Brownian motions (i.e., $\alpha_1=\alpha_2=1$) with $r(s,t)=r$ a constant has been well studied in the literature; see e.g., \cite{Iyengar85} and \cite{Metzller10}. Therein \De{an explicit expression} for  \eqref{eq:Pru} was given through the modified Bessel function and in the form of series; \De{recently} \cite{ShaoWang13} obtained some computable bounds for  \eqref{eq:Pru}. \eP{We refer to \cite{LieMan07} for related results.}\\
Explicit calculation of \eqref{eq:Pru} is only possible for correlated Brownian motions.
\De{Since typically in applications calculation of the joint survival probability is needed for large thresholds $u$, one can rely on the asymptotic theory to  find adequate approximations of this survival probability.} In \cite{PiterStam05} logarithmic asymptotics of  \eqref{eq:Pru}  as $u\rw\IF$ for \De{general correlated} Gaussian processes $X_1, X_2$ was obtained; see also \cite{Debicki10} for a \Ee{general treatment of the multidimensional case.}  \Ee{So far in the literature} there are only two contributions that derive exact asymptotics of \eqref{eq:Pru} \Ee{for certain Gaussian processes}, namely Cheng and Xiao \cite{ChengXiao13} obtained an exact asymptotic expansion of \eqref{eq:Pru} for two correlated smooth Gaussian processes  $X_1, X_2$. In the aforementioned paper the result  was obtained  by studying the geometric properties of the processes.
The second contribution is from Anshin \cite{Anshin05} \Ee{where} the exact asymptotics for two correlated non-smooth Gaussian processes $X_1, X_2$ is \Ee{derived} by relying on a modified double-sum method (see \cite{HRF, PicandsB, PicandsA, Pit72, Pit96} for details on  the double-sum method). The assumptions in \cite{Anshin05} are \Pe{such} that our model of two \De{standard} fBm's with $r(s,t)=r\in(-1,1)$ is not included.
Indeed, \Ee{the conditions} {\bf C1-C3} therein are all invalid \Ee{for our} model. Due to wide applications of fBm's and their exit probabilities, we consider in this paper the exact asymptotics of $P_r(u)$ given as in \eqref{eq:Pru}
with $X_1,X_2$ being the two \De{standard} fBm's above with a constant cross-correlation function $r\in(-1,1)$. \eP{Another merit of choosing fBm's rather than other (general) Gaussian processes is that it allows for somewhat explicit formulae.} In order to proceed with our analysis using a modified double-sum method, we \Ee{shall extend the celebrated} Borell-TIS inequality and \De{the} Piterbarg inequality for 2-dimensional Gaussian random fields in \netheo{Borell} and \netheo{ThmPiter}, respectively. These results are of independent interest \Ee{given their} importance in the theory of Gaussian processes and random fields; \De{see} \cite{MWX} for new developments in this direction.

Before presenting our main \eLP{results} we recall \Ee{the definition of the} well-known Pickands constant
\BQN\label{LaT}
  \mathcal{H}_{\alpha}:=\lim_{T \rightarrow\infty} \frac{1}{T}\H^0_\alpha\De{[\LamT]},\ \ \alpha\in(0,2], \quad \LamT:=[0,T],
\EQN
where
\BQN\label{defH}
\H^b_\alpha[\Lambda]:= \E{ \exp\biggl(\sup_{t\in\Lambda}\Bigl(\sqrt{2}B_{\alpha}(t)-\abs{t}^{\alpha}-b t\Bigr)\biggr)},\ \ \Lambda\subset \R,\ b\inr.
\EQN
Here  $\{B_\alpha(t),t\inr \}$ is a \eP{standard} fBm defined on $\R$ with Hurst index
$\alpha/2 \in (0,1]$.
By the symmetry about 0 of the fBm,
for any $T\in(0,\IF)$ we have (with $\LamT$ given as in \eqref{LaT})
\BQN\label{defH}
\H_\alpha^0[[-T,0]]=\H_\alpha^0[\LamT],\ \ \H^{-b}_1[[-T,0]]= \H^b_1[\LamT],\ \ \ b\inr. 
\EQN
\Ee{We refer to \cite{albin2010new, bermansojourns, debicki2002ruin, debicki2008note, DRolski, Kosi, DikerY, Man07, Pit96} for the basic properties of \Pe{the} Pickands and related constants.}

Our \eLP{first} principle result is presented below.
\BT\label{thmmain}
Let $\{X_i(t),t\ge0\}$, $i=1,2$ be two \eLP{standard} fBm's with Hurst indexes $\alpha_i/2\in(0,1),i=1,2,$ respectively. If \Pe{$(X_1,X_2)$ are jointly Gaussian with 
a constant cross correlation function} $r\in(-1,1)$,  then as $u\rw\IF$
\BQN
\De{P_r(u)}=\frac{(1+r)^{\frac{3}{2}}}{2\pi\sqrt{1-r}}\Upsilon_1(u) \Upsilon_2(u)u^{-2}\exp\LT(-\frac{u^2}{1+r}\RT)(1+o(1)),
\EQN
where
\BQNY
\Upsilon_i(u)=\left\{
            \begin{array}{ll}
2^{1-\frac{1}{\alpha_i}}(1+r)^{1-\frac{2}{\alpha_i}} \frac{1}{\alpha_i} \H_{\alpha_i}u^{\frac{2}{\alpha_i}-2}, & \hbox{if } \alpha_i\in(0,1),\\
\frac{2+r}{1+r}, &  \hbox{if } \alpha_i=1,\\
1, &  \hbox{if } \alpha_i\in(1,2),
               \end{array}
            \right.\ \ i=1,2.
\EQNY
\ET

{\it {\bf Remarks:}
a) The case $r=0$ can be confirmed by the fact that (cf. \cite{DebickiRol02} or \cite{Pit96}), for a standard fBm \Pe{$B_\alpha$}
\BQN\label{eq:fBm}
\pk{\sup_{t\in[0,1]}B_\alpha(t)>u}= \frac{1}{\sqrt{2\pi}} \mathcal{F}_\alpha(u)u^{-1}\exp\LT(-\frac{u^2}{2}\RT)\ooo,\ \ \asu,
\EQN
where $ \mathcal{F}_\alpha(u)$ is equal to $2^{1-1/\alpha}\alpha^{-1}\H_\alpha u^{2/\alpha-2}$ if $\alpha\in(0,1)$, 2 if $\alpha=1$, and 1 if $\alpha\in(1,2)$.

b)
It follows from \netheo{thmmain} and \eqref{eq:fBm} that $\De{M_1}:=\sup_{s\in[0,1]}X_1(s)$ and $\De{M_2}:=\sup_{t\in[0,1]}X_2(t)$ are \Ee{asymptotically} independent, i.e.,
\BQN\label{eq:ShaoWang}
\lim_{u\rw\IF}\frac{\pk{\De{M_1>u,M_2}>u}}{\pk{M_2>u}}=0.
\EQN
We refer to  \cite{Denuitetal05} for the concept of the asymptotic independence of two random variables. The result in \eqref{eq:ShaoWang} improves \Pe{that} in Corollary 1.1 in \cite{ShaoWang13} where an upper bound (1/2) for the left-hand side of \eqref{eq:ShaoWang} for two correlated Brownian motions was obtained.

c) In \netheo{thmmain} we considered the \eP{joint} extremes of two \De{standard} fBm's on \Pe{the} time interval $[0,1]$. \Pe{Next}, we briefly discuss the case where the time interval is $[0,S]$, \eP{with $S$ some} positive constant. It follows by the self-similarity of the fBm's that, for $S^{\alpha_1-\alpha_2}\ge 1$ \Pe{we have}
\BQNY
\pk{\De{\sup}_{t\in[0,S]}X_1(s)>u,\ \De{\sup}_{t\in[0,S]}X_2(t)>u}=\pk{M_1>c(S^{-\frac{\alpha_2}{2}}u),\ \De{M_2}>S^{-\frac{\alpha_2}{2}}u},
\EQNY
with $ c={S^{ \frac{\alpha_2-\alpha_1}{2}}}\in(0,1]$. By slight modifications of the proofs of \netheo{thmmain} and Lemma A we conclude that similar results can be 
obtained for the case $c>r$. However, the case $c\le r$ can not be dealt with similarly since we do not  observe a similar result as  Lemma A which is crucial for the double-sum method. It turns out that the case  $c\le r$ may not be easily solved in general; new techniques will be explored for it elsewhere.}\\

In the framework of ruin theory, see, e.g.,  \cite{AsmAlb10, DHJ13a, EKM97, HP08, HZ08}, \Pe{\cite{DHJ13a}} 
given that ruin happens one wants to know when it happens. With this motivation, we are interested to know when the first passages occur given that $X_1, X_2$ both \eP{ever pass} a threshold $u>0$ on $[0,1]$. 

Define the first passage times of $X_1,X_2$ to the threshold $u$ by
\BQN\label{deftau}
\tau_1(u)=\inf\{s\ge0, X_1(s)>u\}\ \ \ \text{and}\ \ \ \tau_2(u)=\inf\{t\ge0, X_2(t)>u\},
\EQN
respectively (here we use the \Ke{common} assumption that $\inf\{\emptyset\}=\IF$).  Further,
define $\tAu,\  \tAuu, u>0$ in the same probability space
such that
\BQN\label{eq:taustar}
(\tAu, \tAuu)& \equaldis & (\tau_1(u),\tau_2(u))   \Bigl\lvert (\tau_1(u)\le 1, \tau_2(u)\le 1),
\EQN
where $\equaldis$ stands for equality of distribution functions. With motivation from the aforementioned contributions, our second principle result is concerned with the \Ke{distributional} approximation of the random vector  $(\tAu,
\tAuu)$, as $u \to \IF$. Let $E_i, i=1,2$  be two independent unit exponential random variables, and  denote by $\todis$ the convergence in distribution.

\BT\label{thmmain2}
Under the assumptions of \netheo{thmmain} we have as $u\to\IF$
\BQN
\LT(u^2(1-\tAu), u^2(1-\tAuu)\RT)&  \todis& \LT(\frac{2(1+r)}{\alpha_1}E_1, \frac{2(1+r)}{\alpha_2}E_2\RT).
\EQN
\ET

\eLP{
{\it
{\bf Remark:} Let $\De{M_1}, \De{M_2}$ be given as in Remark b) above. By the self-similarity of the fBm's we have for any $x_1,x_2\ge0, u>0$
\BQNY
\pk{M_1>u+\frac{x_1}{u},M_2>u+\frac{x_2}{u} \Bigl\lvert M_1>u,M_2>u}=\frac{\pk{\sup_{s\in[0,1]}X_1(S_{1,u}s)>u,\sup_{t\in[0,1]}X_2(S_{2,u}t)>u}}{\pk{M_1>u,M_2>u}},
\EQNY
where $S_{i,u}=(1+x_i u^{-2})^{-2/\alpha_i}, i=1,2, u>0$. Therefore, by a similar argument as in the proof of \netheo{thmmain2} we conclude that
\BQN\label{eq:MM}
\limit{u}\pk{M_1>u+\frac{x_1}{u},M_2>u+\frac{x_2}{u} \Bigl\lvert M_1>u,M_2>u}&=& \exp\LT(-\frac{x_1+x_2}{1+r}\RT). 
\EQN
\Ke{In view of Theorem 4.1 in \cite{HASH09} (see also Section 4.1 in \cite{HASH12})}
\BQNY
\limit{u}\pk{X_1(1)>u+\frac{x_1}{u},X_2(1)>u+\frac{x_2}{u} \Bigl\lvert X_1(1)>u,X_2(1)>u}&\Ke{=}& \exp\LT(-\frac{x_1+x_2}{1+r}\RT)
\EQNY
\Ke{holds for any $x_1,x_2 \in [0,\IF)$}, from which we see that
\eqref{eq:MM} is not surprising since the minimum of the function $h(s,t), (s,t)\in(0,1]^2$ \eP{given in \eqref{eq:hh}} is attained at the unique point $(1,1)$, at which the processes usually contribute most to the asymptotics.}}\\

Organization of the rest of the paper: In Section 2 we present some preliminary results including the Borell-TIS inequality and the Piterbarg inequality for 2-dimensional Gaussian random fields. \eLP{The proofs of Theorems \ref{thmmain} and \ref{thmmain2} are} given in Section 3, while proofs of other results are relegated to Appendix.

\section{\Ee{Preliminaries}}
\def\Kk{\Ee{\mathcal{K}}}
\def\simS{\De{\tau_m^2}}

In the asymptotic theory of Gaussian processes, two of the important inequalities are the Borell-TIS inequality (cf. \cite{AdlerTaylor, Pit96}) and the Piterbarg inequality (cf. \cite{Pit96}).
Let $\{Z(t), t\in \Kk \}$ be a centered Gaussian process with  a.s. continuous sample paths, and \Pe{let} 
$\Kk\subset\R$ \Pe{be} a compact set \Ee{with Lebesgue} measure mes($\Kk)>0$. The Borell-TIS inequality, which was proved by \cite{Borell75} and \cite{TIS76} independently, states that
\BQN\label{eq:Bo1}
\pk{\sup_{t\in \Kk }Z(t)>u}\le \exp\LT(-\frac{(u-\mu)^2}{2} \simS \RT)
\EQN
holds for any $u\ge \mu:=\E{\sup_{t\in \Kk}Z(t)}$, with $ \simS :=\inf_{t\in \Kk}\LT(\Var{Z(t)}\RT)^{-1}\in(0,\IF)$.\\
The upper bound in \eqref{eq:Bo1} might not be precise enough for \Ee{various applications} \Ee{due to the appearance of the} constant $\mu$.
V.I. Piterbarg obtained an upper bound under a global H\"{o}lder condition on the Gaussian process, which eliminates the constant $\mu$;
see e.g., Theorem 8.1 in \cite{Pit96} or Theorem 8.1 in \cite{Pit2001}. Specifically, \Hh{if there are some positive constants $\gamma$ and $G$} such that
$\E{(Z(t)-Z(t'))^2}\le G \abs{t-t'}^\gamma$ for all $t,t'\in \Kk$, then
\BQN\label{eq:Pit}
\pk{\sup_{t\in \Kk}Z(t)>u}\le C\ \text{mes}(\Kk)\ u^{\frac{2}{\gamma}-1} \exp\LT(-\frac{u^2}{2} \simS \RT)
\EQN
holds for any $u$ large enough, with some positive constant $C$ \Ee{not depending on $u$}. The last inequality is commonly
referred to as the Piterbarg inequality; see e.g., Proposition 3.2 in \cite{MR3062433} for the case of \cLP{chi-processes}.

\def\Vv{\mathcal{V}}
\Ee{Next, let $ \Vv \subset\R^2$} be a compact set, and \Pe{let $\{(Z_1(t), Z_2(t)) , t\ge0\}$ be a 2-dimensional centered vector Gaussian process 
with components which have} a.s.\ continuous sample paths. Motivated by the \Ee{findings} of \cite{Debicki10, PiterStam05}, we present in \netheo{Borell} and \netheo{ThmPiter}  generalizations of the Borell-TIS and Piterbarg inequalities for 2-dimensional Gaussian random fields $\{(Z_1(s),Z_2(t)),(s,t)\in \Vv\}$. As it will be seen from the proof of \netheo{thmmain}, the generalized Borell-TIS and Piterbarg inequalities are very powerful tools. 


\BT\label{Borell} Let $\{Z_i(t), t\ge0\}, i=1,2$ be two centered Gaussian processes with a.s.\ continuous sample paths, variance functions $\sigma_i(t),i=1,2$ \Pe{being further 
jointly Gaussian with cross-correlation} function $r(s,t)\in(-1,1)$.
Then there exists a constant $\mu$ such that for $u\ge\mu$
\BQN\label{eq:Borell}
\pk{\underset{(s,t)\in \Vv} \bigcup \LT\{Z_1(s)>u, Z_2(t)>u\RT\}}\le   \exp\LT(-\frac{(u-\mu)^2}{2} \simS \RT),
\EQN
where 
\COM{\BQNY
 \simS =\inf_{(s,t)\in \Vv} \LT(\frac{1}{\inf(\sigma_1^2(s),\sigma_2^2(t))} +\frac{ \sigma_1^2(s)+\sigma_2^2(t)-2\sigma_1(s)\sigma_2(t)r(s,t) }{\sigma_1^2(s)\sigma_2^2(t)(1-r^2(s,t))}I(r(s,t)<c(s,t)) \RT),\ \text{with}\ 
c(s,t)=\min\LT(\frac{\sigma_1(s)}{\sigma_2(t)},\frac{\sigma_2(t)}{\sigma_1(s)}\RT).
\EQNY
}
$ \simS =\inf_{(s,t)\in \Vv}\sigma^2(s,t)\Ee{>0}$ with \Pe{(below $I(\cdot)$ stands for the indicator function)}
\BQN\label{eq:sigM}
\sigma^2(s,t)=\frac{1}{\min(\sigma_1^2(s),\sigma_2^2(t))}\LT(1+\frac{(c(s,t)-r(s,t))^2}{1-r^2(s,t)}I(r(s,t)<c(s,t))\RT),\ 
c(s,t)=\min\LT(\frac{\sigma_1(s)}{\sigma_2(t)},\frac{\sigma_2(t)}{\sigma_1(s)}\RT) .
\EQN
In particular, if $r(s,t)<c(s,t)$ for all $(s,t)\in \Vv$, then \eqref{eq:Borell} holds with
\BQN\label{eq:sigM2}
 \simS =\inf_{(s,t)\in \Vv} \frac{ \sigma_1^2(s)+\sigma_2^2(t)-2\sigma_1(s)\sigma_2(t)r(s,t) }{\sigma_1^2(s)\sigma_2^2(t)(1-r^2(s,t))}
\EQN
and further, if $r(s,t)\ge c(s,t)$ for all $(s,t)\in \Vv$, then \eqref{eq:Borell} holds with
\BQNY
 \simS =\inf_{(s,t)\in \Vv}  \frac{1}{\min(\sigma_1^2(s),\sigma_2^2(t))}.
\EQNY
\ET

\BT\label{ThmPiter}
Let $\{Z_i(t), t\ge0\}, i=1,2$ be as in \netheo{Borell}. 
Assume that $\sigma_1(s), \sigma_2(t), r(s,t), (s,t)\in \Vv$ are all \cLP{twice} continuously differentiable with respect to their arguments. If there exist some positive constants \Ee{$\gamma$ and $L$} such that the following global H\"{o}lder condition
\BQN\label{eq:holder}
\E{(Z_i(v_i)-Z_i(w_i))^2}\le L \abs{v_i-w_i}^\gamma, \ \ \ i=1,2
\EQN
holds for all $(v_1,v_2), (w_1,w_2)\in \Vv$, then for all $u$ large
\BQN\label{eq:P}
\pk{\underset{(s,t)\in \Vv} \bigcup \LT\{Z_1(s)>u, Z_2(t)>u\RT\}}\le  C \ \text{\mbox{mes}}(\Vv)\  u^{\frac{4}{\gamma}-1} \exp\LT(-\frac{u^2}{2} \simS \RT),
\EQN
where $ \simS $ is given as in \netheo{Borell}, and $C$ is some positive constant \Ee{not depending on $u$}.
\ET

\begin{remark}
 Assume that $\mathcal{G}=\{(s,t)\in \Vv: (s,t)=\text{arg}\inf\sigma(s,t)\}$ is a finite set. Define  $\mathcal{G}_\vn= \underset{(s,t)\in \mathcal{G}} \bigcup ([s-\vn,s+\vn]\times[t-\vn,t+\vn]\cap\cLP{\Vv})$ for any small positive $\vn$. In view of the proof of  Theorem 8.1 in \cite{Pit96}, the claim \eqref{eq:P} still holds if \eqref{eq:holder} is valid for all $(v_1,v_2), (w_1,w_2)\in \mathcal{G}_\vn$ for some small positive $\vn$.
\end{remark}

\Ee{Now,} we come back to our principle problem of finding the exact asymptotics of $P_r(u)$ as $u\rw\IF$. \Ee{In view of the findings} of \cite{Anshin05, Debicki10, PiterStam05} 
we deduce that the constant $\simS$ given in \eqref{eq:sigM2} (restricted to fBm's case) should play \Ee{a crucial} role in the exact asymptotics of $P_r(u)$. Thus, we need to analyze the following function
\BQN\label{eq:hh}
h(s,t)=\frac{t^{\alpha_2}+s^{\alpha_1}-2rs^{\frac{\alpha_1}{2}}t^{\frac{\alpha_2}{2}}}{s^{\alpha_1}t^{\alpha_2}(1-r^2)},\ \ \ s,t\in (0,1].
\EQN
The function $h(s,t), s,t\in(0,1]$ attains its minimum at \Ee{the} unique point $(s_0,t_0)=(1,1)$ and further
$h(1,1)=\frac{2}{1+r}$. \\

\COM{
\prooflem{Lem:maxima} It is easy to check that
$$
\frac{\partial h}{\partial s}(s,t)= \frac{\partial h}{\partial t}(s,t)=0\ \Leftrightarrow \ s^{\alpha_1}=c^2t^{\alpha_2}.
$$
Therefore, the minimum of $h(s,t)$ on the set $(0,1]^2$ is taken  on the curve $L_1=\{(s,t)\in(0,1]^2:  s^{\alpha_1}=c^2t^{\alpha_2}\}$, or on the lines $L_2=\{(s,t)\in(0,1]^2: s=1\}$ or $L_3=\{(s,t)\in(0,1]^2: t=1\}$. Furthermore, we can derive that
on the line $L_1$, the minimum of $h(s,t)$ is taken at $(c^{2/\alpha_1},1)$ with $h(c^{2/\alpha_1},1)=\frac{2}{1+r}$; on the line $L_2$, the minimum of $h(s,t)$ is taken at $(1,1)$ with $h(1,1)=\frac{1+c^2-2cr}{1-r^2}$; on the line $L_3$, the minimum of $h(s,t)$ is taken at $(1,1)$ with $h(1,1)=\frac{1+c^2-2cr}{1-r^2}$ if $c>r$, and at $\left(\left(\frac{c}{r}\right)^{2/\alpha_1},1\right)$ with $h\left(\left(\frac{c}{r}\right)^{2/\alpha_1},1\right)=1$ if $c\le r$. The claim follows by comparing the values on these points. \QED
\bigskip
}

Let $(\hs,\htt):=(\hat{s}_0(u),\hat{t}_0(u)), u>0$ be a family of points in $[0,1]^2$ satisfying $1-\hs\le (\ln u)^2/u^2$ and $1-\htt\le (\ln u)^2/u^2$.
For the use of the double-sum method, we need to deal with the asymptotics of \cLP{the following joint survival function}
\BQNY
\De{R_{\Lambda_1,\Lambda_2}}(u):=\pk{\underset{(s,t)\in K_u} \bigcup \{X_1(s)>u,X_2(t)>u\}},\ \ \ \asu,
\EQNY
where $K_u=(\hs,\htt)+(u^{-2/\alpha_1}\Lambda_1, u^{-2/\alpha_2}\Lambda_2)$ with $\Lambda_i, i=1,2$ two compact sets in $\R$. Here in our notation, for any $\Lambda\in \R$ $a\Lambda:=\{ax: x\in\Lambda\}$, and  for any $\Lambda\in \R^2$ $(x_1,x_2)+\Lambda:=\{(x_1,x_2)+(y_1,y_2): (y_1,y_2)\in\Lambda\}$.

The following lemma can be seen as a generalization of Pickands and Piterbarg lemmas (cf. \cite{MR3091101, pickands1969asymptotic, Pit72, Pit96}) for 2-dimensional Gaussian random fields. \Ee{Its proof} is presented in Appendix.

{\bf Lemma A}. {\it
Let $\{X_i(t),t\ge0\}, i=1,2$ be two \De{standard} fBm's with Hurst indexes $\alpha_i/2\in(0,1/2],i=1,2,$ respectively. Assume further that the joint correlation function  of them is a constant $r\in(-1,1)$. Then as $u\rw\IF$
\BQN\label{eq:lemmaA}
\De{R_{\Lambda_1,\Lambda_2}}(u)=\piter_{\alpha_1}[\Lambda_1]\piter_{\alpha_2}[\Lambda_2]\frac{(1+r)^{\frac{3}{2}}}{2\pi\sqrt{1-r}}u^{-2}\exp\LT(-\frac{u^2}{2}h(\hs,\htt)\RT)(1+o(1)),
\EQN
where $h(\cdot,\cdot)$ is given as in \eqref{eq:hh}, and
\BQNY
\piter_{\alpha_i}[\Lambda_i]=\left\{
            \begin{array}{ll}
\H_{\alpha_i}^0\LT[\Lambda_i\LT(\frac{1}{\sqrt{2}(1+r) }\RT)^{\frac{2}{\alpha_i}}\RT], & \hbox{if } \alpha_i\in(0,1),\\
\H_{1}^{-(1+r)}\LT[\Lambda_i\LT(\frac{1}{\sqrt{2}(1+r) }\RT)^{2}\RT], &  \hbox{if } \alpha_i=1,\\
               \end{array}
            \right.\ \ i=1,2.
\EQNY
}
\COM{\piter_{\alpha_2}[\Lambda_2]=\left\{
            \begin{array}{ll}
\H_{\alpha_2}^0\LT[\Lambda_2\LT(\frac{1}{\sqrt{2}(1+r) }\RT)^{\frac{2}{\alpha_2}}\RT], & \hbox{if } \alpha_2<1,\\
\H_{1}^{-(1+r)}\LT[\Lambda_2\LT(\frac{1}{\sqrt{2}(1+r) }\RT)^{2}\RT], &  \hbox{if }\alpha_2=1.\\
               \end{array}
            \right.
}

\def\kn{k=1,\cdots,n}
\def\lm{l=1,\cdots,m}

\section{Proofs of \eLP{Theorems \ref{thmmain} and \ref{thmmain2}}}

In this section, we \eLP{first present} the proof of \netheo{thmmain} which is based on a tailored double-sum method as in \cite{Anshin05}; see the classical monograph \cite{Pit96} \Ee{for a deep explanation} on the double-sum method. \eLP{Then we present the proof of \netheo{thmmain2}.}

\prooftheo{thmmain} Let $\delta(u)=(\ln u)^2/u^{2}$, and set $D_u=\{(s,t)\in[0,1]^2: 1-s\le \delta(u), 1-t \le \delta(u)\}$. With these notation we have
\BQNY
P_{1,r}(u):=\pk{\underset{(s,t)\in D_u}  \bigcup  \{X_1(s)>u, X_2(t)>u\}}&\le & P_{r}(u)\\
&\le &P_{1,r}(u)+\pk{\underset{(s,t)\in [0,1]^2/D_u} \bigcup  \{X_1(s)>u, X_2(t)>u\}}\\
&=:&P_{1,r}(u)+P_{2,r}(u).
\EQNY
Next, we shall \De{derive} the exact asymptotics of $P_{1,r}(u)$ as $u\rw\IF$, and show that 
\BQN\label{eq:P2r}
P_{2,r}(u)=o(P_{1,r}(u)),\ \ \ \ \ u\rw\IF
\EQN
implying thus
 \BQNY
P_{r}(u)=P_{1,r}(u)(1+o(1))\ \ \ \ \ u\rw\IF.
\EQNY
\De{Next, we derive an} upper bound for $P_{2,r}(u)$ \De{by utilising} the generalized Borell-TIS and Piterbarg inequalities. 
Choose some small $\vn\in(0,1)$  such that
\BQN\label{eq:hatc}
\hat{c}(s,t):=\min\LT(\frac{t^{\alpha_2/2}}{s^{\alpha_1/2}}, \frac{s^{\alpha_1/2}}{t^{\alpha_2/2}}\RT)>r,\ \ \forall (s,t)\in[1-\vn,1]^2.
\EQN
Clearly, for any $u$ positive
\BQNY
P_{2,r}(u)&\le&\pk{\underset{(s,t)\in [0,1]^2/[1-\vn,1]^2} \bigcup \{ X_1(s)>u, X_2(t)>u\}}+\pk{\underset{(s,t)\in [1-\vn,1]^2/D_u} \bigcup \{ X_1(s)>u, X_2(t)>u\}}.
\EQNY
 It follows from the Borell-TIS inequality in \netheo{Borell} that for all $u$ large
\BQN\label{eq:borell}
\pk{\underset{(s,t)\in [0,1]^2/[1-\vn,1]^2} \bigcup \{ X_1(s)>u, X_2(t)>u\}}\le  \exp\LT(-\frac{(u-\mu)^2}{2}\inf_{(s,t)\in(0,1]^2/[1-\vn,1]^2}f(s,t)\RT),
\EQN
where $\mu\in(0,\IF)$ is some constant and
$$
f(s,t)=\frac{1}{\min\LT(s^{\alpha_1}, t^{\alpha_2}\RT)}\LT(1+\frac{(\hatc(s,t)-r)^2}{1-r^2}I(r<\hatc(s,t))\RT),\ \ (s,t)\in (0,1]^2/[1-\vn,1]^2.
$$
Further, straightforward calculations yield that (recall \eqref{eq:hh} for the expression of  $h(\cdot,\cdot)$)
\BQNY
\inf_{(s,t)\in(0,1]^2/[1-\vn,1]^2}f(s,t)> h(1,1)=\frac{2}{1+r}.
\EQNY
Moreover, in view of \eqref{eq:hatc} we have from the Piterbarg inequality in \netheo{ThmPiter} and its remark that,  for all $u$ large
\BQNY
\pk{\underset{(s,t)\in [1-\vn,1]^2/D_u} \bigcup \{ X_1(s)>u, X_2(t)>u\}}\le C u^{\frac{4}{\min(\alpha_1,\alpha_2)}-1}\exp\LT(-\frac{u^2}{2}\inf_{(s,t)\in [1-\vn,1]^2/D_u}h(s,t)\RT),
\EQNY
\Ee{with $C>0$ not depending on $u$.} In addition from the Taylor expansion of $h(s,t)$ around the point $(1,1)$ we have
$$
h(s,t)=h(1,1)+\frac{1}{1+r}\LT(\alpha_1(1-s)+\alpha_2(1-t)\RT)(1+o(1)).
$$
Hence, for the chosen small enough $\vn\Pe{>0}$ there exists some positive constant $C_1$ such that
\BQNY
h(s,t)\ge h(1,1)+C_1 \delta(u)
\EQNY
for any $(s,t)\in [1-\vn,1]^2/D_u$, implying thus,   for all $u$ large
\BQN\label{eq:piter}
\pk{\underset{(s,t)\in [1-\vn,1]^2/D_u} \bigcup \{ X_1(s)>u, X_2(t)>u\}}\le C u^{\frac{4}{\min(\alpha_1,\alpha_2)}-1}\psi_r(u)\exp\LT(-\frac{C_1}{2}(\ln u)^2\RT),
\EQN
\Ee{where we set
$$ \psi_r(u): = \exp\LT(-\frac{u^2}{2} h(1,1)\RT) = \exp\LT(-\frac{u^2}{1+r} \RT) .$$}
Consequently, from \eqref{eq:borell} and \eqref{eq:piter} we obtain the following upper bound for $P_{2,r}(u)$ when $u$ is large
\BQN\label{eq:P2rupper}
P_{2,r}(u)\le \exp\LT(-\frac{(u-\mu)^2}{2}\inf_{(s,t)\in[0,1]^2/[1-\vn,1]^2}f(s,t)\RT)+ C u^{\frac{4}{\min(\alpha_1,\alpha_2)}-1}\psi_r(u)\exp\LT(-\frac{C_1}{2}(\ln u)^2\RT).
\EQN
From now on we focus on the asymptotics of $P_{1,r}(u)$ as $u\rw\IF$.
Let $T_1,T_2$ be two positive constants. For $\alpha_i\le 1, i=1,2$, we can split the rectangle $D_u$ into several sub-rectangles of side lengths $T_1u^{-2/\alpha_1}$ and $T_2u^{-2/\alpha_2}$. Specifically, let
$$
\Del_{k,l}=\Delo_k\times\Delt_l=[s_{k+1},s_k]\times[t_{l+1},t_l],\ \ \ k,l\inn \bigcup \{0\},
$$
with $s_k=1-kT_1u^{-2/\alpha_1}$ and $t_l=1-lT_2u^{-2/\alpha_2}$, and further set
\BQNY
N_i(u)=\LT\lfloor T_i^{-1}(\ln u)^2 u^{\frac{2}{\alpha_i}-2}\RT\rfloor+1,\ \ \ i=1,2.
\EQNY
Here $\lfloor \cdot\rfloor$ denotes the ceiling function. Thus
\BQN
 \bigcup _{k=0}^{N_1(u)-1} \bigcup _{l=0}^{N_2(u)-1} \Del_{k,l} \subset D_u \subset  \bigcup _{k=0}^{N_1(u)} \bigcup _{l=0}^{N_2(u)} \Del_{k,l}.
\EQN
In what follows, we deal with \Pe{only three cases} (distinguished by $\alpha_i$'s):

\underline{Case i) $\alpha_1\in(0,1)$ and $\alpha_2\in(0,1)$}.
Applying the Bonferroni  inequality in Lemma B (given in Appendix) we obtain
\BQNY
P_{1,r}(u)\le \sum_{k=0}^{N_1(u)}\sum_{l=0}^{N_2(u)}\pk{\sup_{s\in\Delo_k}X_1(s)>u,\sup_{t\in\Delt_l}X_2(t)>u}
\EQNY
and
\BQN\label{eq:p1lower}
P_{1,r}(u)&\ge& \sum_{k=0}^{N_1(u)-1}\sum_{l=0}^{N_2(u)-1}\pk{\sup_{s\in\Delo_k}X_1(s)>u,\sup_{t\in\Delt_l}X_2(t)>u}-\Sigma_1(u)-\Sigma_2(u),
\EQN
where
\BQNY
&&\Sigma_1(u)=\sum_{k=0}^{N_1(u)}\underset{0\le l_1< l_2\le N_2(u)}{\sum\sum}\pk{\sup_{s\in\Delo_k}X_1(s)>u,\sup_{t\in\Delt_{l_1}}X_2(t)>u,\sup_{t\in\Delt_{l_2}}X_2(t)>u},\\
&&\Sigma_2(u)=\sum_{l=0}^{N_2(u)}\underset{0\le k_1< k_2\le N_1(u)}{\sum\sum}\pk{\sup_{s\in\Delo_{k_1}}X_1(s)>u,\sup_{s\in\Delo_{k_2}}X_1(s)>u,\sup_{t\in\Delt_{l}}X_2(t)>u}.
\EQNY
Further, in view of Lemma A 
\BQNY
P_{1,r}(u)&\le& \H^0_{\alpha_1}\LT[[-T_1,0]\LT(\frac{1}{\sqrt{2}(1+r)}\RT)^{\frac{2}{\alpha_1}}\RT] \H^0_{\alpha_2}\LT[[-T_2,0]\LT(\frac{1}{\sqrt{2}(1+r)}\RT)^{\frac{2}{\alpha_2}}\RT]\\
&& \times \frac{(1+r)^{\frac{3}{2}}}{2\pi\sqrt{1-r}}u^{-2}\sum_{k=0}^{N_1(u)}\sum_{l=0}^{N_2(u)}\exp\LT(-\frac{u^2}{2}h(s_k,t_l)\RT)(1+o(1))
\EQNY
as $u\rw\IF$. Since by Taylor expansion
$$
h(s_k,t_l)=h(1,1)+\frac{1}{1+r}\LT(\alpha_1(1-s_k)+\alpha_2(1-t_l)\RT)(1+o(1)),\ \ \ u\rw\IF
$$
 we have
\BQNY
\sum_{k=0}^{N_1(u)}\sum_{l=0}^{N_2(u)}\exp\LT(-\frac{u^2}{2}h(s_k, t_l)\RT)&=&
\psi_r(u)\frac{u^{\frac{2}{\alpha_1}+\frac{2}{\alpha_2}-4}}{T_1T_2}
\prod_{j=1}^2\biggl(\int_0^\IF\exp\LT(-\frac{\alpha_j }{2(1+r)} x\RT)dx\biggr)\ooo.
\EQNY
Therefore, as $u\to \IF$
\BQN\label{eq:upperP1r}
P_{1,r}(u)&\le& \frac{2^{1-\frac{1}{\alpha_1}-\frac{1}{\alpha_2}}(1+r)^{\frac{7}{2}-\frac{2}{\alpha_1}-\frac{2}{\alpha_2}}}{\pi \alpha_1\alpha_2\sqrt{1-r}}
\frac{\H^0_{\alpha_1}[0,b_1 T_1]}{b_1 T_1}\frac{\H^0_{\alpha_2}[0,b_2 T_2]}{b_2 T_2}
u^{ \frac{2}{\alpha_1}+\frac{2}{\alpha_2}-6 }
\psi_r(u)(1+o(1)),
\EQN
where $b_i=\LT(1/(\sqrt{2}(1+r))\RT)^{2/\alpha_i}, i=1,2.$
The same arguments yield that
\BQN\label{eq:lowerP1r}
\sum_{k=0}^{N_1(u)-1}\sum_{l=0}^{N_2(u)-1}\pk{\sup_{s\in\Delo_k}X_1(s)>u,\sup_{t\in\Delt_l}X_2(t)>u}&=&\frac{2^{1-\frac{1}{\alpha_1}-\frac{1}{\alpha_2}}(1+r)^{\frac{7}{2}-\frac{2}{\alpha_1}-\frac{2}{\alpha_2}}}{\pi \alpha_1\alpha_2\sqrt{1-r}}
\frac{\H^0_{\alpha_1}[0,b_1 T_1]}{b_1 T_1}\frac{\H^0_{\alpha_2}[0,b_2 T_2]}{b_2 T_2}\nonumber\\
&&\times u^{ \frac{2}{\alpha_1}+\frac{2}{\alpha_2}-6 }
\psi_r(u)(1+o(1))
\EQN
 as $u\rw\IF$. Next, we consider the estimates of $\Sigma_i(u), i=1,2$. To this end, we
define, for any $T,T_0\in(0,\IF)$
 $$
 \H^0_\alpha([0,T],[T_0,T_0+T])=\int_{-\IF}^\IF \exp(x)\pk{\sup_{t\in[0,T]}\sqrt{2}B_{\alpha}(t)-\abs{t}^{\alpha}>x,\sup_{t\in[T_0,T_0+T]}\sqrt{2}B_{\alpha}(t)-\abs{t}^{\alpha}>x }dx,\ \ \alpha\in(0,2)
 $$
and denote, for any $n\ge 1$
 $$
 \H^0_\alpha(n;T)=\H^0_\alpha([0,T],[nT,(n+1)T]).
 $$
 It follows from Lemma 3 in \cite{Anshin05} or Lemmas 6 and 7 in \cite{HusLadPit10} that
 \BQN
 \lim_{T\rw\IF}\frac{\sum_{n=1}^\IF\H^0_\alpha(n;T)}{T}=0.
 \EQN
Since
\BQNY
&&\pk{\sup_{s\in\Delo_k}X_1(s)>u,\sup_{t\in\Delt_{l_1}}X_2(t)>u,\sup_{t\in\Delt_{l_2}}X_2(t)>u}\\
&&=\pk{\sup_{s\in\Delo_k}X_1(s)>u,\sup_{t\in\Delt_{l_1}}X_2(t)>u}
+\pk{\sup_{s\in\Delo_k}X_1(s)>u,\sup_{t\in\Delt_{l_2}}X_2(t)>u}\\
&&-\pk{\sup_{s\in\Delo_k}X_1(s)>u,\sup_{t\in\Delt_{l_1} \bigcup \Delt_{l_2}}X_2(t)>u}
\EQNY
similar arguments as \Pe{in} the derivation of \eqref{eq:upperP1r} \De{imply} that
\BQN\label{eq:Sig1}
\Sigma_1(u)&\le& \frac{2^{1-\frac{1}{\alpha_1}-\frac{1}{\alpha_2}}(1+r)^{\frac{7}{2}-\frac{2}{\alpha_1}-\frac{2}{\alpha_2}}}{\pi \alpha_1\alpha_2\sqrt{1-r}}
\frac{\H^0_{\alpha_1}[0,b_1 T_1]}{b_1 T_1}\sum_{n=1}^\IF\frac{\H^0_{\alpha_2}[n; b_2 T_2]}{b_2 T_2}\nonumber\\
&&\times u^{ \frac{2}{\alpha_1}+\frac{2}{\alpha_2}-6 }
\psi_r(u)(1+o(1)).
\EQN
Similarly
\BQN\label{eq:Sig2}
\Sigma_2(u)&\le& \frac{2^{1-\frac{1}{\alpha_1}-\frac{1}{\alpha_2}}(1+r)^{\frac{7}{2}-\frac{2}{\alpha_1}-\frac{2}{\alpha_2}}}{\pi \alpha_1\alpha_2\sqrt{1-r}}
\sum_{n=1}^\IF\frac{\H^0_{\alpha_1}[n;b_1 T_1]}{b_1 T_1}\frac{\H^0_{\alpha_2}[0,b_2 T_2]}{b_2 T_2}\nonumber\\
&&\times u^{ \frac{2}{\alpha_1}+\frac{2}{\alpha_2}-6 }
\psi_r(u)(1+o(1)).
\EQN
Consequently, from (\ref{eq:upperP1r}-\ref{eq:Sig2}) by letting $T_1,T_2\rw\IF$ \De{we obtain}
\BQN
P_{1,r}(u)&=& \frac{2^{1-\frac{1}{\alpha_1}-\frac{1}{\alpha_2}}(1+r)^{\frac{7}{2}-\frac{2}{\alpha_1}-\frac{2}{\alpha_2}}}{\pi \alpha_1\alpha_2\sqrt{1-r}}
\H_{\alpha_1}\H_{\alpha_2}
u^{ \frac{2}{\alpha_1}+\frac{2}{\alpha_2}-6 }
\psi_r(u)(1+o(1))\ \ \ \asu.
\EQN

\underline{Case ii) $\alpha_1\in(0,1)$ and $\alpha_2=1$}.
Applying the Bonferroni inequality we \De{have}
\BQNY
P_{1,r}(u)&\le& \sum_{k=0}^{N_1(u)} \pk{\sup_{s\in\Delo_k}X_1(s)>u,\sup_{t\in\Delt_0}X_2(t)>u}\\
&&+\sum_{k=0}^{N_1(u)}\sum_{l=1}^{N_2(u)}\pk{\sup_{s\in\Delo_k}X_1(s)>u,\sup_{t\in\Delt_l}X_2(t)>u}
\EQNY
and
\BQNY
P_{1,r}(u)&\ge& \sum_{k=0}^{N_1(u)-1} \pk{\sup_{s\in\Delo_k}X_1(s)>u,\sup_{t\in\Delt_0}X_2(t)>u}-\Sigma_3(u),
\EQNY
where
\BQNY
\Sigma_3(u)=\underset{0\le k_1< k_2\le N_1(u)}{\sum\sum}\pk{\sup_{s\in\Delo_{k_1}}X_1(s)>u,\sup_{s\in\Delo_{k_2}}X_1(s)>u,\sup_{t\in\Delt_{0}}X_2(t)>u}.
\EQNY
By Lemma A 
\BQN\label{eq:upperCase2}
&&\sum_{k=0}^{N_1(u)(\text{or}\ N_1(u)-1)} \pk{\sup_{s\in\Delo_k}X_1(s)>u,\sup_{t\in\Delt_0}X_2(t)>u}\nonumber\\
&&= \H^0_{\alpha_1}[0,b_1T_1] \H^{1+r}_{1}[0,b_2T_2] \frac{(1+r)^{\frac{3}{2}}}{2\pi\sqrt{1-r}}u^{-2}\sum_{k=0}^{N_1(u)}\exp\LT(-\frac{u^2}{2}h(s_k,1)\RT)(1+o(1))\nonumber\\
&&=\frac{2^{-\frac{1}{\alpha_1}}(1+r)^{\frac{5}{2}-\frac{2}{\alpha_1}}}{\pi \alpha_1\sqrt{1-r}}
\frac{\H^0_{\alpha_1}[0,b_1 T_1]}{b_1 T_1}\H^{1+r}_{1}[0,b_2 T_2]u^{ \frac{2}{\alpha_1}-4 }
\psi_r(u)(1+o(1))
\EQN
as $u\rw\IF$, where $b_i, i=1,2$ are the same as in \eqref{eq:upperP1r}. Similarly
\BQN\label{eq:lowerCase2}
&&\sum_{k=0}^{N_1(u)}\sum_{l=1}^{N_2(u)}\pk{\sup_{s\in\Delo_k}X_1(s)>u,\sup_{t\in\Delt_l}X_2(t)>u}\nonumber\\
&&=\frac{2^{-\frac{1}{\alpha_1}}(1+r)^{\frac{5}{2}-\frac{2}{\alpha_1}}}{\pi \alpha_1\sqrt{1-r}}
\frac{\H^0_{\alpha_1}[0,b_1 T_1]}{b_1 T_1}\H^{0}_{1}[0,b_2 T_2]\sum_{l=1}^{\IF}\exp\LT(-\frac{ T_2 l}{2(1+r)}\RT)u^{ \frac{2}{\alpha_1}-4 }
\psi_r(u)(1+o(1))
\EQN
as $u\rw\IF$. Moreover, it follows with similar arguments as in \eqref{eq:Sig1} that
\BQN\label{eq:Sig3}
\Sigma_3(u)\le\frac{2^{-\frac{1}{\alpha_1}}(1+r)^{\frac{5}{2}-\frac{2}{\alpha_1}}}{\pi \alpha_1\sqrt{1-r}}
\sum_{n=1}^\IF\frac{\H^0_{\alpha_1}[n;b_1 T_1]}{b_1 T_1}\H^{1+r}_{1}[0,b_2 T_2]u^{ \frac{2}{\alpha_1}-4 }
\psi_r(u)(1+o(1))
\EQN
as $u\rw\IF$. Consequently, letting $T_1,T_2\rw\IF$ from (\ref{eq:upperCase2}-\ref{eq:Sig3}) we have
\BQNY
P_{1,r}(u)=\frac{2^{-\frac{1}{\alpha_1}}(2+r)(1+r)^{\frac{3}{2}-\frac{2}{\alpha_1}}}{\pi \alpha_1\sqrt{1-r}}
\H_{\alpha_1} u^{ \frac{2}{\alpha_1}-4 }
\psi_r(u)(1+o(1))\ \ \ \asu,
\EQNY
where we used the fact that $\H_1^{1+r}\cLP{:=\lim_{T\to\IF}\H_1^{1+r}[\Lambda_T]}=(2+r)/(1+r)$; see e.g., \cite{DeMan03} or \cite{MR3091101}.

\underline{Case iii) $\alpha_1\in(0,1)$ and $\alpha_2\in(1,2)$}.
Since $\alpha_2>1,$ it follows that $\delta(u)\subset\Delt_0$. Thus
\BQNY
P_{1,r}(u)&\le& \sum_{k=0}^{N_1(u)} \pk{\sup_{s\in\Delo_k}X_1(s)>u,\sup_{t\in\Delt_0}X_2(t)>u}
\EQNY
and further
\BQNY
P_{1,r}(u)&\ge& \sum_{k=0}^{N_1(u)-1} \pk{\sup_{s\in\Delo_k}X_1(s)>u, X_2(1)>u}-\Sigma_4(u),
\EQNY where
\BQNY
\Sigma_4(u)= \underset{0\le k_1< k_2\le N_1(u)}{\sum\sum}\pk{\sup_{s\in\Delo_{k_1}}X_1(s)>u,\sup_{s\in\Delo_{k_2}}X_1(s)>u,X_2(1)>u}.
\EQNY
Using the same technique as in the proof of Lemma A 
 (or let $T_2\to 0$ therein), we can show that
\BQNY
\pk{\sup_{s\in\Delo_k}X_1(s)>u, X_2(1)>u}&=&  \H^0_{\alpha_1}\LT[[-T_1,0]\LT(\frac{1}{\sqrt{2}(1+r)}\RT)^{\frac{2}{\alpha_1}}\RT] \\ 
&& \times \frac{(1+r)^{\frac{3}{2}}}{2\pi\sqrt{1-r}}u^{-2} \exp\LT(-\frac{u^2}{2}h(s_k,1)\RT)(1+o(1))
\EQNY
as $u\rw\IF$, implying
\BQN\label{eq:upperCase3}
\sum_{k=0}^{N_1(u)-1} \pk{\sup_{s\in\Delo_k}X_1(s)>u, X_2(1)>u} =\frac{2^{-\frac{1}{\alpha_1}}(1+r)^{\frac{5}{2}-\frac{2}{\alpha_1}}}{\pi \alpha_1\sqrt{1-r}}
\frac{\H^0_{\alpha_1}[0,b_1 T_1]}{b_1 T_1} u^{ \frac{2}{\alpha_1}-4 }
\psi_r(u)(1+o(1)), \quad u\to \IF.
\EQN
Moreover
\BQN\label{eq:Sig4}
\Sigma_4(u)\le\frac{2^{-\frac{1}{\alpha_1}}(1+r)^{\frac{5}{2}-\frac{2}{\alpha_1}}}{\pi \alpha_1\sqrt{1-r}}
\sum_{n=1}^\IF\frac{\H^0_{\alpha_1}[n;b_1 T_1]}{b_1 T_1} u^{ \frac{2}{\alpha_1}-4 }
\psi_r(u)(1+o(1))
\EQN
as $u\rw\IF$.  Consequently, letting $T_1\rw\IF, T_2\rw0$ we conclude from \eqref{eq:upperCase2}, \eqref{eq:upperCase3} and \eqref{eq:Sig4} that
\BQNY
P_{1,r}(u) =\frac{2^{-\frac{1}{\alpha_1}}(1+r)^{\frac{5}{2}-\frac{2}{\alpha_1}}}{\pi \alpha_1\sqrt{1-r}}
 \H_{\alpha_1}  u^{ \frac{2}{\alpha_1}-4 }
 \psi_r(u)(1+o(1)) \ \ \ \asu.
\EQNY
With all the techniques used in the proofs of Cases i)-iii) we see that the other cases  \Pe{for the possible choices of $\alpha_1$ and $\alpha_2$} can be shown similarly without any further difficulty, thus the detailed proofs are omitted. Moreover, it follows from \eqref{eq:P2rupper} and the asymptotics of $P_{1,r}(u)$ in any of the remaining cases 
that \eqref{eq:P2r} holds, and thus the proof is complete. \QED

\eLP{\prooftheo{thmmain2}  First note that, for any $x_1,x_2\ge0, u>0$
\BQNY
\pk{u^2(1-\tAu)>x_1, u^2(1-\tAuu)>x_2}=\frac{\pk{\sup_{s\in[0,T_{1,u}]}X_1(s)>u,\ \sup_{t\in[0,T_{2,u}]}X_2(t)>u}}{\pk{\sup_{s\in[0,1]}X_1(s)>u,\ \sup_{t\in[0,1]}X_2(t)>u}},
\EQNY
with $T_{i,u}=1-x_i u^{-2},  i=1,2$. Further, we \eP{write}
\BQNY
\pk{\sup_{s\in[0,T_{1,u}]}X_1(s)>u,\ \sup_{t\in[0,T_{2,u}]}X_2(t)>u}&=&\pk{\sup_{s\in[0,1]}\widetilde{X_1}(s)>u,\ \sup_{t\in[0,1]}\widetilde{X_2}(t)>u},
\EQNY
where $\widetilde{X_i}(t):=X_i(T_{i,u}t), t\in[0,1]$. Define $\widetilde{h_u}(s,t):=h(T_{1,u}s,T_{2,u}t), (s,t)\in(0,1]^2$, with $h(\cdot,\cdot)$ given as in \eqref{eq:hh}. It follows from a slight modification of the proof of Lemma A  that \eqref{eq:lemmaA} holds for $\widetilde{X_1}, \widetilde{X_2}$, without any other changes apart from that $h(\cdot,\cdot)$ is replaced by $\widetilde{h_u}(\cdot,\cdot)$. With this modification of Lemma A, by a similar argument as in the proof of \netheo{thmmain}  we conclude that, as $u\to\IF$
\BQN
\pk{\sup_{s\in[0,1]}\widetilde{X_1}(s)>u,\ \sup_{t\in[0,1]}\widetilde{X_2}(t)>u}
=\frac{(1+r)^{\frac{3}{2}}}{2\pi\sqrt{1-r}}\Upsilon_1(u) \Upsilon_2(u)u^{-2}\exp\LT(-\frac{u^2}{2} \widetilde{h_u}(1,1)\RT)(1+o(1)),
\EQN
where $\Upsilon_i(u), i=1,2$ are given \eP{as} in \netheo{thmmain}. Consequently, from the last formula and \netheo{thmmain}, for any $x_1,x_2\ge0$
\BQNY
\pk{u^2(1-\tAu)>x_1, u^2(1-\tAuu)>x_2}&=&\exp\LT(-\frac{u^2}{2} (\widetilde{h_u}(1,1)-h(1,1))\RT)\oo\\
&\to& \exp\LT(-\LT(\frac{\alpha_1}{2(1+r)}x_1+\frac{\alpha_2}{2(1+r)}x_2 \RT)\RT),\ \ \ u\to\IF
\EQNY
establishing thus the claim, and hence the proof is complete.
\QED
}

\section{Appendix}
\Hh{Below we present} the proofs of \netheo{Borell}, \netheo{ThmPiter} and Lemma A. \Hh{
We also state and prove Lemma B which is of some interest on its own.}

\prooftheo{Borell} Denote
$$
A(s,t)=\sigma_1^2(s)+\sigma_2^2(t)-2\sigma_1(s)\sigma_2(t)r(s,t),\ \ \ (s,t)\in \Vv.
$$
Next, we introduce two nonnegative functions $a(s,t), b(s,t),  (s,t)\in \Vv$ \Pe{as follows}
\BQNY
a(s,t)=\left\{
            \begin{array}{ll}
\frac{\sigma_2^2(t)-\sigma_1(s)\sigma_2(t)r(s,t) }{A(s,t)} , & \hbox{if } c(s,t)>r(s,t), 
\\
1, & \hbox{if } c(s,t)\le r(s,t) \ \text{and}\  \sigma_1(s)\le \sigma_2(t),\\
0,&  \hbox{otherwise},
              \end{array}
            \right.  \ \ \  (s,t)\in \Vv
\EQNY and
\BQNY
b(s,t)=\left\{
            \begin{array}{ll}
\frac{\sigma_1^2(s)-\sigma_1(s)\sigma_2(t)r(s,t) }{A(s,t)} , & \hbox{if } c(s,t)>r(s,t),\\ 
1,& \hbox{if } c(s,t)\le r(s,t) \ \text{and}\  \sigma_2(t)< \sigma_1(s),\\
0,&  \hbox{otherwise},
              \end{array}
            \right.\ \ \  (s,t)\in \Vv.
\EQNY
Since $a(s,t)+b(s,t)\cLP{=1},  (s,t)\in \Vv$, it follows that
\BQN\label{eq:BorellY}
\lefteqn{\pk{\underset{(s,t)\in \Vv} \bigcup \Bigl\{Z_1(s)>u, Z_2(t)>u\Bigr\}}}\nonumber\\
&\le& \pk{\underset{(s,t)\in \Vv} \bigcup \{ a(s,t)Z_1(s)+b(s,t)Z_2(t)>a(s,t) u+ b(s,t)u\}}\nonumber
\\
&=& \pk{\sup_{(s,t)\in \Vv} Y(s,t;a,b) >u}
\EQN
where
$$ Y(s,t;a,b)= a(s,t)Z_1(s)+b(s,t)Z_2(t),\ \ \ \ (s,t)\in \Vv.$$
\Ee{Since  further}
\BQNY
(\E{(Y(s,t;a,b))^2})^{-1}&=&\frac{1}{a^2(s,t)\sigma_1^2(s)+b^2(s,t)\sigma_2^2(t)+2a(s,t)b(s,t)\sigma_1(s)\sigma_2(t)r(s,t)}\\
&=&\frac{ \sigma_1^2(s)+\sigma_2^2(t)-2\sigma_1(s)\sigma_2(t)r(s,t) }{\sigma_1^2(s)\sigma_2^2(t)(1-r^2(s,t))}I(c(s,t)>r(s,t))\\
&&+ \frac{1}{\min(\sigma_1^2(s),\sigma_2^2(t))} I(c(s,t)\le r(s,t))
\EQNY
the claim follows from the Borell-TIS inequality for one-dimensional Gaussian random fields (e.g., \cite{AdlerTaylor}) with
$$\mu=\E{\sup_{(s,t)\in \Vv} Y(s,t;a,b)}<\IF$$
and thus the proof is complete.  \QED

\prooftheo{ThmPiter} We use  the same notation as in the proof of \netheo{Borell}. In the light of \eqref{eq:BorellY} and Theorem 8.1 in \cite{Pit96}, it suffices to show that
\BQN \label{eq:regu}
\E{(Y(s,t;a,b)-Y(s',t';a,b))^2}\le L_1 (\abs{s-s'}^\gamma+\abs{t-t'}^\gamma),\ \ \forall (s,t), (s',t')\in \Vv
\EQN
holds for some positive constants $L_1$ and $\gamma$, which can be confirmed by some straightforward calculations,  and thus the claim follows. \QED

{\bf Proof of Lemma A}: Using the classical technique, see e.g., \eLP{\cite{Anshin05, HusLadPit10, Pit96}}, we have for any $u>0$
\BQN\label{eq:uxy}
\De{R_{\Lambda_1,\Lambda_2}}(u)&=&\frac{1}{u^2}\int_{-\IF}^\IF\int_{-\IF}^\IF\pk{\underset{(s,t)\in K_u} \bigcup \{X_1(s)>u,X_2(t)>u\}|X_1(\hs)=u-\frac{x}{u},X_2(\htt)=u-\frac{y}{u}}\\
&&\times f_{X_1(\hs),X_2(\htt)}\left(u-\frac{x}{u},u-\frac{y}{u}\right)dxdy, \nonumber
\EQN
where
\BQNY
&& f_{X_1(\hs),X_2(\htt)}\left(u-\frac{x}{u},u-\frac{y}{u}\right)=\frac{1}{2\pi\sqrt{\hs^{\alpha_1}\htt^{\alpha_2}(1-r^2)}}
\exp\LT(-\frac{1}{2\hs^{\alpha_1}\htt^{\alpha_2}(1-r^2) u^2}\LT(\htt^{\alpha_2}x^2+\hs^{\alpha_1}y^2-2r\hs^{\alpha_1/2} \htt^{\alpha_2/2}xy\RT)\RT)\\
&& \times
\exp\LT(-\frac{1}{2\hs^{\alpha_1}\htt^{\alpha_2}(1-r^2)  }\LT(-2\htt^{\alpha_2}x-2\hs^{\alpha_1}y +2r\hs^{\alpha_1/2} \htt^{\alpha_2/2}\LT(x+y\RT)\RT)\RT)\exp
 \LT(-\frac{u^2}{2} h(\hs,\htt)\RT),\ \ x,y\inr,
\EQNY
where $h(\cdot,\cdot)$ is defined as in \eqref{eq:hh}. \Pe{Set} for $x,y\inr$
$$
\xi_u(s)=u(X_1(\hs+u^{-\frac{2}{\alpha_1}}s)-u)+x, \ \ \eta_u(t)=u(X_2(\htt+u^{-\frac{2}{\alpha_2}}t)-u)+y.
$$
The probability in the integrand of \eqref{eq:uxy} can be rewritten as
$$
p_u(x,y)=\pk{\underset{(s,t)\in \Lambda_1\times\Lambda_2} \bigcup \{\xi_u(s)>x,\eta_u(t)>y\}|\xi_u(0)=0,\eta_u(0)=0}.
$$
\Pe{Next, we calculate the expectation and covariance of the conditional random vector} $(\xi_u(s),\eta_u(t))|(\xi_u(0),\eta_u(0))$. We have
\BQNY
\E{\begin{array}{c}
     \xi_u(s) \\
     \eta_u(t)
   \end{array}\Bigg| \begin{array}{c}
     \xi_u(0) \\
     \eta_u(0)
   \end{array}
} = \E{\begin{array}{c}
     \xi_u(s) \\
     \eta_u(t)
   \end{array}
}+A\left(
     \begin{array}{c}
       \xi_u(0)-\E{\xi_u(0)} \\
        \eta_u(0)-\E{\eta_u(0)} \\
     \end{array}
   \right),
\EQNY
where
\BQNY
A=\Cov\LT(\left(
\begin{array}{c}
     \xi_u(s) \\
     \eta_u(t)
   \end{array}
          \right),
          \left(
\begin{array}{c}
     \xi_u(0) \\
     \eta_u(0)
   \end{array}
          \right)
\RT)\times \Cov\LT(
\begin{array}{c}
     \xi_u(0) \\
     \eta_u(0)
   \end{array}
         \RT)^{-1}
\EQNY
and further
\BQNY
\Cov\LT(\begin{array}{c}
          \xi_u(t)-\xi_u(s) \\
          \eta_u(t_1)-\eta_u(s_1)
        \end{array}\Bigg|
        \begin{array}{c}
          \xi_u(0) \\
          \eta_u(0)
        \end{array}
\RT)=\Cov\LT(\begin{array}{c}
          \xi_u(t)-\xi_u(s) \\
          \eta_u(t_1)-\eta_u(s_1)
        \end{array}
\RT)+B\ \Cov\left(
            \begin{array}{c}
              \xi_u(0) \\
              \eta_u(0) \\
            \end{array}
          \right)^{-1} B^\top,
\EQNY
where
$$
B=\left(
    \begin{array}{cc}
      b_{11}(u) & b_{12}(u) \\
      b_{21}(u) & b_{22}(u) \\
    \end{array}
  \right)=
\Cov\LT(\LT(\begin{array}{c}
          \xi_u(t)-\xi_u(s) \\
          \eta_u(t_1)-\eta_u(s_1)
        \end{array}\RT),\LT(
        \begin{array}{c}
          \xi_u(0) \\
          \eta_u(0)
        \end{array}
\RT)\RT).
$$
Further 
\COM{\BQNY
\Cov\LT(\left(
\begin{array}{c}
     \xi_u(s) \\
     \eta_u(t)
   \end{array}
          \right),
          \left(
\begin{array}{c}
     \xi_u(0) \\
     \eta_u(0)
   \end{array}
          \right)
\RT)=u^2\left(
          \begin{array}{cc}
            \frac{1}{2}\LT(\hs^{\alpha_1}+(\hs+u^{-\frac{2}{\alpha_1}}s)^{\alpha_1}-u^{-2}s^{\alpha_1}\RT) & r\htt^{\frac{\alpha_2}{2}} \LT(\hs+u^{-\frac{2}{\alpha_1}}s\RT)^{\frac{\alpha_1}{2}}\\
            r\hs^{\frac{\alpha_1}{2}} \LT(\htt+u^{-\frac{2}{\alpha_2}}t\RT)^{\frac{\alpha_2}{2}} &  \frac{1}{2}\LT(\htt^{\alpha_2}+(\htt+u^{-\frac{2}{\alpha_2}}t)^{\alpha_2}-u^{-2}t^{\alpha_2}\RT) \\
          \end{array}
        \right).
\EQNY
Therefore
}
\BQN\label{eq:Ainv}
\Cov\left(
            \begin{array}{c}
              \xi_u(0) \\
              \eta_u(0) \\
            \end{array}
          \right)^{-1}=\frac{u^{-2}}{\htt^{\alpha_2}\hs^{\alpha_1}(1-r^2)}
\left(
  \begin{array}{cc}
    \htt^{\alpha_2} & -r\htt^{\frac{\alpha_2}{2}}\hs^{\frac{\alpha_1}{2}} \\
    -r\htt^{\frac{\alpha_2}{2}}\hs^{\frac{\alpha_1}{2}} & \hs^{\alpha_1} \\
  \end{array}
\right)
\EQN and
\BQNY
&&A=\frac{1}{\htt^{\alpha_2}\hs^{\alpha_1}(1-r^2)}\times\\
&&\times\left(
          \begin{array}{cc}
            \frac{1}{2}\htt^{\alpha_2}\LT(\hs^{\alpha_1}+(\hs+u^{-\frac{2}{\alpha_1}}s)^{\alpha_1}-u^{-2}s^{\alpha_1}\RT)- & -\frac{1}{2}r\htt^{\frac{\alpha_2}{2}}\hs^{\frac{\alpha_1}{2}}\LT(\hs^{\alpha_1}+(\hs+u^{-\frac{2}{\alpha_1}}s)^{\alpha_1}-u^{-2}s^{\alpha_1}\RT)+\\
             -r^2 \htt^{ \alpha_2 } \hs^{\frac{\alpha_2}{2}}\LT(\hs+u^{-\frac{2}{\alpha_1}}s\RT)^{\frac{\alpha_1}{2}}&
             +r \htt^{ \frac{\alpha_2}{2}} \hs^{\alpha_1 }\LT(\hs+u^{-\frac{2}{\alpha_1}}s\RT)^{\frac{\alpha_1}{2}}\\
              -\frac{1}{2}r\htt^{\frac{\alpha_2}{2}}\hs^{\frac{\alpha_1}{2}}\LT(\htt^{\alpha_2}+(\htt+u^{-\frac{2}{\alpha_2}}t)^{\alpha_2}-u^{-2}t^{\alpha_2}\RT)+&
               \frac{1}{2}\hs^{\alpha_1}\LT(\htt^{\alpha_2}+(\htt+u^{-\frac{2}{\alpha_2}}t)^{\alpha_2}-u^{-2}t^{\alpha_2}\RT)-\\
            +r\hs^{\frac{\alpha_1}{2}} \htt^{\alpha_2}\LT(\htt+u^{\frac{2}{\alpha_2}}t\RT)^{\frac{\alpha_2}{2}} &  -r^2\htt^{\frac{\alpha_2}{2}}\hs^{\alpha_1}\LT(\htt+u^{-\frac{2}{\alpha_2}}t\RT)^{\frac{\alpha_2}{2}} \\
          \end{array}
        \right).
\EQNY
Set next
\BQNY
\left(
    \begin{array}{cc}
      e_{1 }(u)   \\
      e_{2 }(u)  \\
    \end{array}
  \right):=\E{\begin{array}{c}
     \xi_u(s) \\
     \eta_u(t)
   \end{array}\Bigg| \begin{array}{c}
     \xi_u(0)=0 \\
     \eta_u(0)=0
   \end{array}
}=\left(
    \begin{array}{c}
      x-u^2 \\
      y-u^2 \\
    \end{array}
  \right)+A\left(
             \begin{array}{c}
               u^2-x \\
               u^2-y \\
             \end{array}
           \right)
.
\EQNY
It follows that
\BQNY
e_{1 }(u)=\frac{1}{\htt^{\alpha_2}\hs^{\alpha_1}(1-r^2)}\LT(-\LT(\frac{1}{2} \htt^{\alpha_2}-\frac{1}{2}r\htt^{\frac{\alpha_2}{2}}\hs^{\frac{\alpha_1}{2}}\RT)\abs{s}^{\alpha_1}+\lambda_1(u)u^2+\lambda_2(u)x+\lambda_3(u)y\RT),
\EQNY
where
\BQNY
&&\lambda_1(u)=\frac{1}{2}\htt^{\frac{\alpha_2}{2}}\LT(\htt^{\frac{\alpha_2}{2}}\LT((\hs+u^{-\frac{2}{\alpha_1}}s)^{\frac{\alpha_1}{2}}
+\hs^{\frac{\alpha_1}{2}}-2r^2\hs^{\frac{\alpha_1}{2}}\RT)-r\hs^{\frac{\alpha_1}{2}}\LT((\hs+u^{-\frac{2}{\alpha_1}}s)^{\frac{\alpha_1}{2}}
-\hs^{\frac{\alpha_1}{2}}\RT)\RT)\\
&&\qquad\qquad\ \ \times\LT((\hs+u^{-\frac{2}{\alpha_1}}s)^{\frac{\alpha_1}{2}}
-\hs^{\frac{\alpha_1}{2}}\RT)\\
&&\lambda_2(u)=\frac{1}{2}\htt^{\alpha_2}\LT((\hs+u^{-\frac{2}{\alpha_1}}s)^{\frac{\alpha_1}{2}}
+\hs^{\frac{\alpha_1}{2}}-2r^2\hs^{\frac{\alpha_1}{2}}\RT)\LT(\hs^{\frac{\alpha_1}{2}}-(\hs+u^{-\frac{2}{\alpha_1}}s)^{\frac{\alpha_1}{2}}
\RT)\\
&&\lambda_3(u)=\frac{1}{2}r\htt^{\frac{\alpha_2}{2}}\hs^{\frac{\alpha_1}{2}} \LT(\hs^{\frac{\alpha_1}{2}}-(\hs+u^{-\frac{2}{\alpha_1}}s)^{\frac{\alpha_1}{2}}
\RT)^2.
\EQNY
Further
\BQNY
e_{2 }(u)=\frac{1}{\htt^{\alpha_2}\hs^{\alpha_1}(1-r^2)}\LT(-\LT(\frac{1}{2} \hs^{\alpha_1}-\frac{1}{2} r\htt^{\frac{\alpha_2}{2}}\hs^{\frac{\alpha_1}{2}}\RT)\abs{t}^{\alpha_2}+\delta_1(u)u^2+\delta_2(u)x+\delta_3(u)y\RT),
\EQNY
where
\BQNY
&&\delta_1(u)=\frac{1}{2}\hs^{\frac{\alpha_1}{2}}\LT(\hs^{\frac{\alpha_1}{2}}\LT((\htt+u^{-\frac{2}{\alpha_2}}t)^{\frac{\alpha_2}{2}}
+\htt^{\frac{\alpha_2}{2}}-2r^2\htt^{\frac{\alpha_2}{2}}\RT)-r\htt^{\frac{\alpha_2}{2}}\LT((\htt+u^{-\frac{2}{\alpha_2}}t)^{\frac{\alpha_2}{2}}
-\htt^{\frac{\alpha_2}{2}}\RT)\RT)\\
&&\qquad\qquad\ \ \times\LT((\htt+u^{-\frac{2}{\alpha_2}}t)^{\frac{\alpha_2}{2}}
-\htt^{\frac{\alpha_2}{2}}\RT)\\
&&\delta_2(u)=\frac{1}{2} r\htt^{\frac{\alpha_2}{2}}\hs^{\frac{\alpha_1}{2}} \LT(\htt^{\frac{\alpha_2}{2}}-(\htt+u^{-\frac{2}{\alpha_2}}t)^{\frac{\alpha_2}{2}}
\RT)^2\\
&&\delta_3(u)=\frac{1}{2}\hs^{ \alpha_1}\LT((\htt+u^{-\frac{2}{\alpha_2}}t)^{\frac{\alpha_2}{2}}
+\htt^{\frac{\alpha_2}{2}}-2r^2\htt^{\frac{\alpha_2}{2}}\RT)\LT(\htt^{\frac{\alpha_2}{2}}-(\htt+u^{-\frac{2}{\alpha_2}}t)^{\frac{\alpha_2}{2}}
\RT).
\EQNY
Thus we have
\BQN
\lim_{u\rw\IF}e_1(u)=\left\{
            \begin{array}{ll}
-\frac{1}{2(1+r)}\abs{s}^{\alpha_1} , & \hbox{if } \alpha_1\in(0,1),\\
-\frac{1}{2(1+r)}\abs{s}+\frac{1}{2}s&  \hbox{if } \alpha_1=1
              \end{array}
            \right.
\EQN
and
\BQN
\lim_{u\rw\IF}e_2(u)=\left\{
            \begin{array}{ll}
-\frac{1}{2(1+r)}\abs{t}^{\alpha_2} , & \hbox{if } \alpha_2\in(0,1),\\
-\frac{1}{2(1+r)}\abs{t}+\frac{1}{2 }t&  \hbox{if } \alpha_2=1.
              \end{array}
            \right.
\EQN
\COM{
ii) $c\le r$. In this case we have $(\hs,\htt)\rw(\LT(c/r\RT)^{2/\alpha_1},1)$. Thus
\BQN
\lim_{u\rw\IF}e_1(u)=\left\{
            \begin{array}{ll}
0 , & \hbox{if } \alpha_1<1,\\
\frac{r^2}{2}s&  \hbox{if } \alpha_1=1,
              \end{array}
            \right.
\EQN
and
\BQN
\lim_{u\rw\IF}e_2(u)=\left\{
            \begin{array}{ll}
-\frac{1 }{2}\abs{t}^{\alpha_2} , & \hbox{if } \alpha_2<1,\\
-\frac{1 }{2}\abs{t}+ \frac{1 }{2}t&  \hbox{if } \alpha_2=1.
              \end{array}
            \right.
\EQN
}

Similarly
\BQNY
&&b_{11}(u)= \frac{1}{2}\LT(s^{\alpha_1}-t^{\alpha_1}+u^2\LT(\LT(\hs+u^{-\frac{2}{\alpha_1}}t\RT)^{\alpha_1}-\LT(\hs+u^{-\frac{2}{\alpha_1}}s\RT)^{\alpha_1}\RT)\RT),\\
&&b_{12}(u)= u^2r \htt^{\frac{\alpha_2}{2}}\LT(\LT(\hs+u^{-\frac{2}{\alpha_1}}t\RT)^{\frac{\alpha_1}{2}}-\LT(\hs+u^{-\frac{2}{\alpha_1}}s\RT)^{\frac{\alpha_1}{2}}\RT), \\
&&b_{21}(u)=u^2r \hs^{\frac{\alpha_1}{2}}\LT(\LT(\htt+u^{-\frac{2}{\alpha_2}}t_1\RT)^{\frac{\alpha_2}{2}}-\LT(\htt+u^{-\frac{2}{\alpha_2}}s_1\RT)^{\frac{\alpha_2}{2}}\RT),\\
&&b_{22}(u)=\frac{1}{2}\LT(s_1^{\alpha_2}-t_1^{\alpha_2}+u^2\LT(\LT(\htt+u^{-\frac{2}{\alpha_2}}t_1\RT)^{\alpha_2}-\LT(\htt+u^{-\frac{2}{\alpha_2}}s_1\RT)^{\alpha_2}\RT)\RT),
\EQNY
\COM{
\BQNY
B=\left(
    \begin{array}{cc}
      \frac{c^2}{2}\LT(s^{\alpha_1}-t^{\alpha_1}+u^2\LT(\LT(\hs+u^{-\frac{2}{\alpha_1}}t\RT)^{\alpha_1}-\LT(\hs+u^{-\frac{2}{\alpha_1}}s\RT)^{\alpha_1}\RT)\RT) & cu^2r \htt^{\frac{\alpha_2}{2}}\LT(\LT(\hs+u^{-\frac{2}{\alpha_1}}t\RT)^{\frac{\alpha_1}{2}}-\LT(\hs+u^{-\frac{2}{\alpha_1}}s\RT)^{\frac{\alpha_1}{2}}\RT) \\
      cu^2r \hs^{\frac{\alpha_1}{2}}\LT(\LT(\htt+u^{-\frac{2}{\alpha_2}}t_1\RT)^{\frac{\alpha_2}{2}}-\LT(\htt+u^{-\frac{2}{\alpha_2}}s_1\RT)^{\frac{\alpha_2}{2}}\RT) &  \frac{1}{2}\LT(s_1^{\alpha_2}-t_1^{\alpha_2}+u^2\LT(\LT(\htt+u^{-\frac{2}{\alpha_2}}t_1\RT)^{\alpha_2}-\LT(\htt+u^{-\frac{2}{\alpha_2}}s_1\RT)^{\alpha_2}\RT)\RT)\\
    \end{array}
  \right)
\EQNY
}
which together with \eqref{eq:Ainv} gives that
\BQNY
B\ \Cov\left(
            \begin{array}{c}
              \xi_u(0) \\
              \eta_u(0) \\
            \end{array}
          \right)^{-1} B^\top=
\left(
  \begin{array}{cc}
    o(1) & o(1) \\
    o(1) & o(1) \\
  \end{array}
\right)
\EQNY
as $u\rw\IF.$ Further
\BQNY
&&\Cov(\xi_u(t)-\xi_u(s),\xi_u(t)-\xi_u(s))=\abs{t-s}^{\alpha_1},\ \ \Cov(\eta_u(t_1)-\eta_u(s_1),\eta_u(t_1)-\eta_u(s_1))=\abs{t_1-s_1}^{\alpha_2},\\
&&\Cov(\xi_u(t)-\xi_u(s),\eta_u(t_1)-\eta_u(s_1))=u^2r\LT((\hs+u^{-\frac{2}{\alpha_1}}t)^{\frac{\alpha_1}{2}}-
(\hs+u^{-\frac{2}{\alpha_1}}s)^{\frac{\alpha_1}{2}}\RT)\\
&&\qquad \qquad \qquad \times \LT((\htt+u^{-\frac{2}{\alpha_2}}t_1)^{\frac{\alpha_2}{2}}-
(\htt+u^{-\frac{2}{\alpha_2}}s_1)^{\frac{\alpha_2}{2}}\RT)=o(1),
\EQNY
as $u\rw\IF.$ Therefore,
\BQNY
\Cov\LT(\begin{array}{c}
          \xi_u(t)-\xi_u(s) \\
          \eta_u(t_1)-\eta_u(s_1)
        \end{array}\Bigg|
        \begin{array}{c}
          \xi_u(0)=0 \\
          \eta_u(0)=0
        \end{array}
\RT)=\left(
       \begin{array}{cc}
         \abs{t-s}^{\alpha_1} & o(1) \\
         o(1) & \abs{t_1-s_1}^{\alpha_2} \\
       \end{array}
     \right), \ \ \asu.
\EQNY
Consequently, using similar arguments as in \cite{Anshin05} (see also \cite{debicki2002ruin}, \cite{HusLadPit10} or \cite{Pit96}) we obtain
\BQNY
\lim_{u\rw\IF}p_u(x,y)=
\pk{\sup_{s\in \Lambda_1}\chi_1(s)>x}\pk{\sup_{t\in \Lambda_2}\chi_2(t)>y}
\EQNY
for any $x,y\inr$, where $\chi_1$ and $\chi_2$ are two independent stochastic processes given by
\BQNY
\chi_1(s)=\widehat{B}_{\alpha_1}(s)+\left\{
            \begin{array}{ll}
-\frac{1}{2(1+r)}\abs{s}^{\alpha_1} , & \hbox{if }\alpha_1\in(0,1),\\
-\frac{1}{2(1+r)}\abs{s}+\frac{1}{2}s&  \hbox{if } \alpha_1=1,
              \end{array}
            \right.\ \ \ \ \ s\in\R
\EQNY
and
\BQNY
\chi_2(t)=\widetilde{B}_{\alpha_2}(t)+\left\{
            \begin{array}{ll}
-\frac{1}{2(1+r)}\abs{t}^{\alpha_2} , & \hbox{if } \alpha_2\in(0,1),\\
-\frac{1}{2(1+r)}\abs{t}+\frac{1}{2}t&  \hbox{if }  \alpha_2=1
              \end{array}
            \right.\ \ \ \ \ t\in\R.
\EQNY
Here $\widehat{B}_{\alpha_1}$ and $ \widetilde{B}_{\alpha_2}$ are two independent fBm's defined on $\R$ with Hurst indexes $\alpha_1/2$ \Hh{and} $\alpha_2/2\in(0,1)$, respectively.
Similar arguments as in \cite{Anshin05} and \cite{HusLadPit10} show that the limit (letting  $u\rw\IF$) can be passed \Hh{under} the integral sign in \eqref{eq:uxy}.
It follows then that
\BQNY
&\De{R_{\Lambda_1,\Lambda_2}}(u)=& \ooo\frac{1}{2\pi \sqrt{1-r^2} u^2}\exp\LT(-\frac{u^2}{2}h(\hs,\htt)\RT)\nonumber\\
&&\times \prod_{i=1}^2 \Biggl(\int_{-\IF}^\IF\exp\LT(\frac{x}{1+r}\RT)\pk{\sup_{s\in \Lambda_i}\Hh{\chi_i}(s)>x}dx \Biggr),\ \ u\to \IF.
\EQNY
Since
\BQNY
 \int_{-\IF}^\IF\exp\LT(\frac{x}{1+r}\RT)\pk{\sup_{s\in \Lambda_i}\Hh{\chi_i}(s)>x}dx =\left\{
            \begin{array}{ll}
(1+r)\H_{\alpha_i}^0\LT[\Lambda_1\LT(\frac{1}{\sqrt{2}(1+r) }\RT)^{\frac{2}{\alpha_i}}\RT], & \hbox{if } \alpha_i\in(0,1),\\
(1+r)\H_{1}^{-(1+r)}\LT[\Lambda_i\LT(\frac{1}{\sqrt{2}(1+r) }\RT)^{2}\RT], &  \hbox{if } \alpha_i=1\\
               \end{array}
            \right.
\EQNY
the claim follows. \QED
\COM{and
\BQNY
\int_{-\IF}^\IF\exp\LT(\frac{y}{1+r}\RT)\pk{\sup_{t\in \Lambda_2}\chi_2(t)>y}dy =\left\{
            \begin{array}{ll}
(1+r)\H_{\alpha_2}^0\LT[\Lambda_2\LT(\frac{1}{\sqrt{2}(1+r) }\RT)^{\frac{2}{\alpha_2}}\RT], & \hbox{if } \alpha_2<1,\\
(1+r)\H_{1}^{1+r}\LT[\Lambda_2\LT(\frac{1}{\sqrt{2}(1+r) }\RT)^{2}\RT], &  \hbox{if }\alpha_2=1.\\
               \end{array}
            \right.
\EQNY
Consequently, inserting the above formula into \eqref{eq:Ku} completes the proof. \QED
}

{\bf Lemma B}. {\it Let $(\Omega, \mathfrak{F}, \mathbb{P})$ be a probability space and $A_1,\cdots,A_n$ and $B_1,\cdots,B_m$ be $n+m$ events in $\mathfrak{F}$ for $n,m\ge 2$. Then
\BQN
&&\pk{\underset{\lm}{\underset{\kn} \bigcup }(A_k\cap B_l)}\ge\sum_{k=1}^n\sum_{l=1}^m\pk{A_k\cap B_l}\nonumber\\
&&\ \ \ \qquad -\sum_{k=1}^n\sum_{1\le l_1<l_2\le m}\pk{A_k\cap B_{l_1}\cap B_{l_2}}
-\sum_{l=1}^m\sum_{1\le k_1<k_2\le n}\pk{ A_{k_1}\cap A_{k_2}\cap B_l}.
\EQN
}

{\bf Proof of Lemma B}: The proof relies on the following Bonferroni inequality; see e.g., Lemma 2 in \cite{Michna09}.
\BQNY
\sum_{k=1}^n\pk{A_k}\ge\pk{ \bigcup _{k=1}^n A_k}\ge\sum_{k=1}^n\pk{A_k}
-\sum_{1\le k_1<k_2\le n}\pk{ A_{k_1}\cap A_{k_2}}.
\EQNY
Since further
\BQNY
&&\pk{\underset{\lm}{\underset{\kn} \bigcup }(A_k\cap B_l)}=\pk{ \bigcup _{k=1}^n(A_k\cap( \bigcup _{l=1}^m B_l))}\\
&&\ \ \ \ \ge\sum_{k=1}^n\pk{A_k \cap( \bigcup _{l=1}^m B_l)}
-\sum_{1\le k_1<k_2\le n}\pk{ A_{k_1}\cap A_{k_2} \cap( \bigcup _{l=1}^m B_l)}\\
&&\ \ \ \ \ge\sum_{k=1}^n\sum_{l=1}^m\pk{A_k\cap B_l}-\sum_{k=1}^n\sum_{1\le l_1<l_2\le m}\pk{A_k\cap B_{l_1}\cap B_{l_2}}
-\sum_{l=1}^m\sum_{1\le k_1<k_2\le n}\pk{ A_{k_1}\cap A_{k_2}\cap B_l}
\EQNY
the proof is complete. \QED

\bigskip

{\bf Acknowledgement}: We would like to thank both the {Editor, associate Editor  and the Referees} for their kind comments and corrections 
which significantly \Pe{improved} the manuscript.
We are thankful to {Krzysztof D\c{e}bicki}, Yimin Xiao and Yuzhen Zhou for various discussions related to the topic of this paper.
This work was support from the Swiss National Science Foundation Grant 200021-140633/1 and the project RARE -318984 (an FP7 Marie Curie IRSES Fellowship). 

\bibliographystyle{plain}

 \bibliography{gausbibU}

\newcommand{\nosort}[1]{}\def\polhk#1{\setbox0=\hbox{#1}{\ooalign{\hidewidth
  \lower1.5ex\hbox{`}\hidewidth\crcr\unhbox0}}}
  \def\polhk#1{\setbox0=\hbox{#1}{\ooalign{\hidewidth
  \lower1.5ex\hbox{`}\hidewidth\crcr\unhbox0}}} \def\cprime{$'$}
  \def\cprime{$'$} \def\cprime{$'$} \def\cprime{$'$}
\begin{thebibliography}{10}

\bibitem{AdlerTaylor}
R.J. Adler and J.E. Taylor.
\newblock {\em Random fields and geometry}.
\newblock Springer Monographs in Mathematics. Springer, New York, 2007.

\bibitem{albin2010new}
J.M.P. Albin and H.~Choi.
\newblock A new proof of an old result by {P}ickands.
\newblock {\em Electronic Communications in Probability}, 15:339--345, 2010.

\bibitem{Anshin05}
A.B. Anshin.
\newblock On the probability of simultaneous extrema of two {G}aussian
  nonstationary processes.
\newblock {\em Teor. Veroyatn. Primen.}, 50(3):417--432, 2005.

\bibitem{AsmAlb10}
S.~Asmussen and H.~Albrecher.
\newblock {\em Ruin probabilities}.
\newblock Advanced Series on Statistical Science \& Applied Probability, 14.
  World Scientific Publishing Co. Pte. Ltd., Hackensack, NJ, second edition,
  2010.

\bibitem{bermansojourns}
M.S. Berman.
\newblock Sojourns and extremes of stochastic processes.
\newblock {\em Wadsworth \& Brooks/Cole, Boston}, 1992.

\bibitem{Borell75}
C.~Borell.
\newblock The {B}runn-{M}inkowski inequality in {G}auss space.
\newblock {\em Invent. Math.}, 30(2):207--216, 1975.

\bibitem{ChengXiao13}
D.~Cheng and Y.~Xiao.
\newblock Geometry and excursion probability of multivariate {G}aussian random
  fields.
\newblock {\em Manuscript}, 2014.

\bibitem{debicki2008note}
K.~D\c{e}bicki and P.~Kisowski.
\newblock A note on upper estimates for pickands constants.
\newblock {\em Statistics \& Probability Letters}, 78(14):2046--2051, 2008.

\bibitem{DRolski}
K.~D\c{e}bicki, Z.~Michna, and T.~Rolski.
\newblock Simulation of the asymptotic constant in some fluid models.
\newblock {\em Stoch. Models}, 19(3):407--423, 2003.

\bibitem{debicki2002ruin}
K.~D{\c{e}}bicki.
\newblock Ruin probability for {G}aussian integrated processes.
\newblock {\em Stochastic Processes and their Applications}, 98(1):151--174,
  2002.

\bibitem{DHJ13a}
K.~D{\polhk{e}}bicki, E.~Hashorva, and L.~Ji.
\newblock Gaussian risk model with financial constraints.
\newblock {\em Scandinavian Actuarial Journal, in press}, 2014.

\bibitem{Kosi}
K.~D{\c{e}}bicki and K.~Kosi\'{n}ski.
\newblock On the infimum attained by the reflected fractional {B}rownian
  motion.
\newblock Extremes, to appear, 2014.

\bibitem{Debicki10}
K.~D{\polhk{e}}bicki, K.~M. Kosi{\'n}ski, M.~Mandjes, and T.~Rolski.
\newblock Extremes of multidimensional {G}aussian processes.
\newblock {\em Stochastic Process. Appl.}, 120(12):2289--2301, 2010.

\bibitem{DeMan03}
K.~D{\polhk{e}}bicki and M.~Mandjes.
\newblock Exact overflow asymptotics for queues with many {G}aussian inputs.
\newblock {\em J. Appl. Probab.}, 40(3):704--720, 2003.

\bibitem{DebickiRol02}
K.~D{\polhk{e}}bicki and T.~Rolski.
\newblock A note on transient {G}aussian fluid models.
\newblock {\em Queueing Systems, Theory and Applications}, 42:321--342, 2002.

\bibitem{Denuitetal05}
M.~Denuit, J.~Dhaene, M.~Goovaerts, and R.~Kaas.
\newblock {\em Actuarial Theory for Dependent Risks: Measures, Orders and
  Models}.
\newblock John Wiley \& Sons, Ltd, England, 2005.

\bibitem{DikerY}
A.B. Dieker and B.~Yakir.
\newblock On asymptotic constants in the theory of {G}aussian processes.
\newblock {\em Bernoulli, to appear}, 2014.

\bibitem{EKM97}
P.~Embrechts, C.~Kl{\"u}ppelberg, and T.~Mikosch.
\newblock {\em Modelling extremal events}, volume~33 of {\em Applications of
  Mathematics (New York)}.
\newblock Springer-Verlag, Berlin, 1997.

\bibitem{HRF}
M.~Falk, J.~H\"usler, and R.-D. Reiss.
\newblock Laws of small numbers: extremes and rare events.
\newblock In {\em DMV Seminar}, volume~23. Birkh\"auser, Basel, third edition,
  2010.

\bibitem{HASH09}
E.~Hashorva.
\newblock Asymptotics for {K}otz type {III} elliptical distributions.
\newblock {\em Statist. Probab. Lett.}, 79(7):927--935, 2009.

\bibitem{HASH12}
E.~Hashorva.
\newblock Exact tail asymptotics in bivariate scale mixture models.
\newblock {\em Extremes}, 15(1):109--128, 2012.

\bibitem{MR3091101}
E.~Hashorva, L.~Ji, and V.I. Piterbarg.
\newblock On the supremum of {$\gamma$}-reflected processes with fractional
  {B}rownian motion as input.
\newblock {\em Stochastic Process. Appl.}, 123(11):4111--4127, 2013.

\bibitem{HusLadPit10}
J.~H{\"u}sler, A.~Ladneva, and V.I. Piterbarg.
\newblock On clusters of high extremes of {G}aussian stationary processes with
  {$\epsilon$}-separation.
\newblock {\em Electron. J. Probab.}, 15:no. 59, 1825--1862, 2010.

\bibitem{HP08}
J.~H{\"u}sler and V.~I. Piterbarg.
\newblock A limit theorem for the time of ruin in a {G}aussian ruin problem.
\newblock {\em Stochastic Process. Appl.}, 118(11):2014--2021, 2008.

\bibitem{HZ08}
J.~H{\"u}sler and Y.~Zhang.
\newblock On first and last ruin times of {G}aussian processes.
\newblock {\em Statist. Probab. Lett.}, 78(10):1230--1235, 2008.

\bibitem{Iyengar85}
S.~Iyengar.
\newblock Hitting lines with two-dimentional {B}rownian motion.
\newblock {\em SIAM J. Appl. Math.}, 45:983--989, 1985.

\bibitem{LieMan07}
P.~Lieshout and M.~Mandjes.
\newblock Tandem {B}rownian queues.
\newblock {\em Math. Methods Oper. Res.}, 66(2):275--298, 2007.

\bibitem{Man07}
M.~Mandjes.
\newblock {\em Large deviations for {G}aussian queues}.
\newblock John Wiley \& Sons Ltd., Chichester, 2007.

\bibitem{MWX}
M.M. Meerschaert, W.~Wang, and Y.~Xiao.
\newblock Fernique-type inequalities and moduli of continuity for anisotropic
  {G}aussian random fields.
\newblock {\em Trans. Amer. Math. Soc.}, 365:1081--1107, 2013.

\bibitem{Metzller10}
A.~Metzler.
\newblock On the first passage problem for correlated {B}rownian motion.
\newblock {\em Statist. Prob. Lett.}, 80:277--284, 2010.

\bibitem{Michna09}
Z.~Michina.
\newblock Remarks on {P}ickands theorem.
\newblock {\em http://arxiv.org/pdf/0904.3832.pdf}.

\bibitem{PicandsB}
J.~Pickands, III.
\newblock Maxima of stationary {G}aussian processes.
\newblock {\em Z. Wahrscheinlichkeitstheorie und Verw. Gebiete}, 7:190--223,
  1967.

\bibitem{pickands1969asymptotic}
J.~Pickands, III.
\newblock Asymptotic properties of the maximum in a stationary {G}aussian
  process.
\newblock {\em Trans. Amer. Math. Soc}, 145:75--86, 1969.

\bibitem{PicandsA}
J.~Pickands, III.
\newblock Upcrossing probabilities for stationary {G}aussian processes.
\newblock {\em Trans. Amer. Math. Soc.}, 145:51--73, 1969.

\bibitem{Pit72}
V.I. Piterbarg.
\newblock On the paper by {J}. {P}ickands ``{U}pcrossing probabilities for
  stationary {G}aussian processes''.
\newblock {\em Vestnik Moskov. Univ. Ser. I Mat. Meh.}, 27(5):25--30, 1972.

\bibitem{Pit96}
V.I. Piterbarg.
\newblock {\em Asymptotic methods in the theory of {G}aussian processes and
  fields}, volume 148 of {\em Translations of Mathematical Monographs}.
\newblock American Mathematical Society, Providence, RI, 1996.
\newblock Translated from the Russian by V.V. Piterbarg, revised by the author.

\bibitem{Pit2001}
V.I. Piterbarg.
\newblock Large deviations of a storage process with fractional {B}rowanian
  motion as input.
\newblock {\em Extremes}, 4:147--164, 2001.

\bibitem{PiterStam05}
V.I. Piterbarg and B.~Stamatovich.
\newblock Rough asymptotics of the probability of simultaneous high extrema of
  two {G}aussian processes: the dual action functional.
\newblock {\em Uspekhi Mat. Nauk}, 60(1(361)):171--172, 2005.

\bibitem{ShaoWang13}
J.~Shao and X.~Wang.
\newblock Estimates of the exit probability for two correlated {B}rownian
  motions.
\newblock {\em Adv. in Appl. Probab.}, 45:37--50, 2013.

\bibitem{MR3062433}
Z.~Tan and E.~Hashorva.
\newblock Exact asymptotics and limit theorems for supremum of stationary
  {$\chi$}-processes over a random interval.
\newblock {\em Stochastic Process. Appl.}, 123(8):2983--2998, 2013.

\bibitem{TIS76}
B.S. Tsirelson, I.A. Ibragimov, and V.N. Sudakov.
\newblock Norms of {G}aussian sample functions.
\newblock In {\em Proceedings of the {T}hird {J}apan-{USSR} {S}ymposium on
  {P}robability {T}heory ({T}ashkent, 1975)}, pages 20--41. Lecture Notes in
  Math., Vol. 550, Berlin, 1976. Springer.

\end{thebibliography}
\end{document}